\documentclass{jnmp}

\usepackage{amsmath}

\JNMPnumberwithin{equation}{section} \resetfootnoterule

\newtheorem{theorem}{Theorem}
\newtheorem{lemma}{Lemma}

\theoremstyle{definition}
\newtheorem{definition}{Definition}

\begin{document}

\renewcommand{\evenhead}{J H van der Walt}
\renewcommand{\oddhead}{First Order Cauchy Problems}

\thispagestyle{empty}

\FirstPageHead{*}{*}{****}{\pageref{firstpage}--\pageref{lastpage}}{Article}

\copyrightnote{2007}{J H van der Walt}

\Name{Generalized Solutions to Nonlinear First Order Cauchy
Problems}

\label{firstpage}

\Author{Jan Harm van der Walt}

\Address{Department of Mathematics and Applied Mathematics,
University of Pretoria, Pretoria 0002, South Africa \\
E-mail: janharm.vanderwalt@up.ac.za\\[10pt]}

\Date{Received September 20, 2007; Accepted in Revised Form Month
*, ****}
\begin{abstract}
\noindent The recent significant enrichment, see \cite{vdWalt5}
through \cite{vdWalt6}, of the Order Completion Method for
nonlinear systems of PDEs \cite{Obergugenberger and Rosinger}
resulted in the global existence of generalized solutions to a
large class of such equations.  In this paper we consider the
existence and regularity of the generalized solutions to a family
of nonlinear first order Cauchy problems. The spaces of
generalized solutions are obtained as the completion of
suitably constructed uniform convergence spaces.\\
\end{abstract}

\quotation{``...provided also if need be that the notion of a
solutions shall
be suitably extended"}\\
\\
--cited form Hilbert's 20th problem \\

\section{Introduction}

It is widely held misconception that there can be no general, type
independent theory for the existence and regularity of solutions
to nonlinear PDEs.  Arnold \cite{Arnold} ascribes this to the more
complicated geometry of $\mathbb{R}^{n}$, as apposed to
$\mathbb{R}$, which is relevant to ODEs alone.  Evans
\cite{Evans}, on the other hand, cites the wide variety of
physical and probabilistic phenomena that are modelled with PDEs.

There are, however, two general, type independent theories for the
solutions of nonlinear PDEs.  The Central Theory for PDEs, as
developed by Neuberger \cite{Neuberger 1}, is based on a
generalized method of steepest descent in suitably constructed
Hilbert spaces.  It delivers generalized solutions to nonlinear
PDEs in a type independent way, although the method is
\textit{not} universally applicable.  However, it does yield
spectacular numerical results.  The Order Completion Method
\cite{Obergugenberger and Rosinger}, on the other hand, yields the
generalized solutions to arbitrary, continuous nonlinear PDEs of
the form
\begin{eqnarray}
T\left(x,D\right)u\left(x\right)=f\left(x\right)\label{PDE}
\end{eqnarray}
where $\Omega\subseteq\mathbb{R}^{n}$ is open and nonempty, $f$ is
continuous, and the PDE operator $T\left(x,D\right)$ is defined
through some jointly continuous mapping
\begin{eqnarray}
F:\Omega\times\mathbb{R}^{K}\rightarrow \mathbb{R}
\end{eqnarray}
by
\begin{eqnarray}
T\left(x,D\right):u\left(x\right)\mapsto
F\left(x,u\left(x\right),...,D^{\alpha}u\left(x\right),..\right)
\end{eqnarray}
The generalized solutions are obtained as elements of the Dedekind
completion of certain spaces of functions, and may be assimilated
with usual Hausdorff continuous, interval valued functions on the
domain of the PDE operator \cite{Anguelov and Rosinger 1}.

Recently, see \cite{vdWalt5} through \cite{vdWalt7}, the Order
Completion Method \cite{Obergugenberger and Rosinger} was
\textit{reformulated} and \textit{enriched} by introducing
suitable uniform convergence spaces, in the sense of
\cite{Beattie}.  In this new setting it is possible, for instance,
to treat PDEs with addition smoothness, over and above the mere
\textit{continuity} of the PDE operator, in a way that allows for
a significantly higher degree of regularity of the solutions
\cite{vdWalt6}.

The aim of this paper is to show how the ideas developed in
\cite{vdWalt5} through \cite{vdWalt7} may be applied to initial
and / or boundary value problems.  In this regard, we consider a
family of nonlinear first order Cauchy problems.  The generalized
solutions are obtained as elements of the completion of a suitably
constructed uniform convergence space.  We note the relative ease
and simplicity of the method presented here, compared to the usual
linear function analytic methods. In this way we come to note
another of the advantages in solving initial and/or boundary value
problems for linear and nonlinear PDEs in this way. Namely,
initial and/or boundary value problems are solved by precisely the
same kind of constructions as the free problems.  On the other
hand, as is well known, this is not so when function analytic
methods - in particular, involving distributions, their
restrictions to lower dimensional manifolds, or the associated
trace operators - are used for the solution of such problems.

The paper is organized as follows.  In Section 2 we introduce some
definitions and results as are required in what follows.  We omit
the proofs, which can be found in \cite{vdWalt5} through
\cite{vdWalt7}.  Section 3 is concerned with the solutions of a
class of nonlinear first order Cauchy problems.  In Section 4 we
discuss the possible interpretation of the generalized solutions
obtained.

\section{Preliminaries}

Let $\Omega$ be some open and nonempty subset of $\mathbb{R}^{n}$,
and let $\overline{\mathbb{R}}$ denote the extended real line
\begin{eqnarray}
\overline{\mathbb{R}}=\mathbb{R}\cup\{\pm\infty\}\nonumber
\end{eqnarray}
A function $u:\Omega\rightarrow \overline{\mathbb{R}}$ belongs to
$\mathcal{ML}^{m}_{0}\left(\Omega\right)$, for some integer $m$,
whenever $u$ is normal lower semi-continuous, in the sense of
Dilworths \cite{Dilworth}, and
\begin{eqnarray}
\begin{array}{ll}
\exists & \Gamma_{u}\subset\Omega\mbox{ closed nowhere dense :} \\
& \begin{array}{ll}
1) & \textit{mes}\left(\Gamma_{u}\right)=0 \\
2) & u\in\mathcal{C}^{m}\left(\Omega\setminus\Gamma_{u}\right) \\
\end{array}\\
\end{array}\label{ML0Def}
\end{eqnarray}
Here $\textit{mes}\left(\Gamma_{u}\right)$ denotes the Lebegue
measure of the set $\Gamma_{u}$.  Recall \cite{Anguelov} that a
function $u:\Omega\rightarrow \overline{\mathbb{R}}$ is normal
lower semi-continuous whenever
\begin{eqnarray}
\begin{array}{ll}
\forall & x\in\Omega\mbox{ :} \\
& I\left(S\left(u\right)\right)\left(x\right)=u\left(x\right) \\
\end{array}\label{NLSCDef}
\end{eqnarray}
where
\begin{eqnarray}
\begin{array}{ll}
\forall & u:\Omega\rightarrow \overline{\mathbb{R}}\mbox{ :} \\
& \begin{array}{ll} 1) & I\left(u\right):\Omega\ni x\mapsto
\sup\{\inf\{u\left(y\right)\mbox{ : }y\in
B_{\delta}\left(x\right)\}\mbox{ : }\delta>0\}\in\mathbb{\overline{R}} \\
2) & S\left(u\right):\Omega\ni x\mapsto
\inf\{\sup\{u\left(y\right)\mbox{ : }y\in
B_{\delta}\left(x\right)\}\mbox{ : }\delta>0\}\in\mathbb{\overline{R}} \\
\end{array} \\
\end{array}\nonumber
\end{eqnarray}
are the lower- and upper- Baire Operators, respectively, see
\cite{Anguelov} and \cite{Baire}.  Note that each function
$u\in\mathcal{ML}^{m}_{0}\left(\Omega\right)$ is measurable and
nearly finite with respect to Lebesgue measure.  In particular,
the space $\mathcal{ML}^{m}_{0}\left(\Omega\right)$ contains
$\mathcal{C}^{m}\left(\Omega\right)$.  In this regard, we note
that the partial differential operators
\begin{eqnarray}
D^{\alpha}:\mathcal{C}^{m}\left(\Omega\right)\rightarrow
\mathcal{C}^{0}\left(\Omega\right)\mbox{, }|\alpha|\leq m\nonumber
\end{eqnarray}
extend to mappings
\begin{eqnarray}
\mathcal{D}^{\alpha}:\mathcal{ML}^{m}_{0}\left(\Omega\right)\ni u
\rightarrow \left(I\circ S\right)\left(D^{\alpha}u\right)\in
\mathcal{ML}^{0}_{0}\left(\Omega\right)\label{PDOps}
\end{eqnarray}

A convergence structure $\lambda_{a}$, in the sense of
\cite{Beattie}, may be defined on
$\mathcal{ML}^{0}_{0}\left(\Omega\right)$ as follows.
\begin{definition}\label{CAEDef}
For any $u\in \mathcal{ML}^{0}_{0}\left(\Omega\right)$, and any
filter $\mathcal{F}$ on $\mathcal{ML}^{0}_{0}\left(\Omega\right)$,
\begin{eqnarray}
\mathcal{F}\in\lambda_{a}\left(u\right)\Leftrightarrow
\left(\begin{array}{ll}
\exists & E\subset\Omega\mbox{ :} \\
& \begin{array}{ll}
a) & \textit{mes}\left(E\right)=0 \\
b) & x\in\Omega\setminus E \Rightarrow \mathcal{F}\left(x\right)\mbox{ converges to }u\left(x\right) \\
\end{array} \\
\end{array}\right)\nonumber
\end{eqnarray}
Here $\mathcal{F}\left(x\right)$ denotes the filter of real
numbers given by
\begin{eqnarray}
\mathcal{F}\left(x\right)=[\{\{v\left(x\right)\mbox{ : }v\in
F\}\mbox{ : }F\in\mathcal{F}\}]
\end{eqnarray}
\end{definition}
That $\lambda_{a}$ does in fact constitute a uniform convergence
structure on $\mathcal{ML}^{0}_{0}\left(\Omega\right)$ follows by
\cite[Example 1.1.2 (iii)]{Beattie}.  Indeed, $\lambda_{a}$ is the
almost everywhere convergence structure, which is Hausdorff. One
may now introduce a \textit{complete} uniform convergence
structure $\mathcal{J}_{a}$, in the sense of \cite{Beattie}, on
$\mathcal{ML}^{0}_{0}\left(\Omega\right)$ in such a way that the
induced convergence structure \cite[Definition 2.1.3]{Beattie} is
$\lambda_{a}$.
\begin{definition}\label{UCAEDef}
A filter $\mathcal{U}$ on $\mathcal{ML}^{0}_{0}\left(\Omega\right)
\times \mathcal{ML}^{0}_{0}\left(\Omega\right)$ belongs to
$\mathcal{J}_{a}$ whenever there exists $k\in\mathbb{N}$ such that
\begin{eqnarray}
\begin{array}{ll}
\forall & i=1,...,k\mbox{ :} \\
\exists & u_{i}\in\mathcal{ML}^{0}_{0}\left(\Omega\right)\mbox{ :} \\
\exists & \mathcal{F}_{i}\mbox{ a filter on $\mathcal{ML}^{0}_{0}\left(\Omega\right)$ :} \\
& \begin{array}{ll}
a) & \mathcal{F}_{i}\in\lambda_{a}\left(u_{i}\right) \\
b) & \left(\mathcal{F}_{1}\times\mathcal{F}_{1}\right)\cap...\cap\left(\mathcal{F}_{k}\times\mathcal{F}_{k}\right)\subseteq\mathcal{U} \\
\end{array} \\
\end{array}\nonumber
\end{eqnarray}
\end{definition}
The uniform convergence structure $\mathcal{J}_{a}$ is referred to
as the uniform convergence structure associated with the
convergence structure $\lambda_{a}$, see \cite[Proposition
2.1.7]{Beattie}.

We note that the concept of a convergence structure on a set $X$
is a generalization of that of topology on $X$.  With every
topology $\tau$ on $X$ one may associate a convergence structure
$\lambda_{\tau}$ on $X$ through
\begin{eqnarray}
\begin{array}{ll}
\forall & x\in X\mbox{ :} \\
\forall & \mathcal{F}\mbox{ a filter on $X$ :} \\
& \mathcal{F}\in\lambda_{\tau}\left(x\right) \Leftrightarrow \mathcal{V}_{\tau}\left(x\right)\subseteq\mathcal{F} \\
\end{array}\nonumber
\end{eqnarray}
where $\mathcal{V}_{\tau}\left(x\right)$ denotes the filter of
$\tau$-neighborhoods at $x$.  However, not every convergence
structure $\lambda$ on $X$ is induced by a topology in this way.
Indeed, the convergence structure $\lambda_{a}$ specified above is
one such an example.  A uniform convergence space is the
generalization of a uniform space in the context of convergence
spaces.  The reader is referred to \cite{Beattie} for details
concerning convergence spaces.

\section{First Order Cauchy Problems}

Let $\Omega=\left(-a,a\right)\times
\left(-b,b\right)\subset\mathbb{R}^{2}$, for some $a,b>0$, be the
domain of the independent variables $\left(x,y\right)$.  We are
given
\begin{eqnarray}
F:\overline{\Omega}\times\mathbb{R}^{4}\rightarrow \mathbb{R}
\label{DefFunc}
\end{eqnarray}
and
\begin{eqnarray}
f:[-a,a]\rightarrow \mathbb{R}\label{ICFunc}
\end{eqnarray}
Here $F$ is jointly continuous in all of its variables, and $f$ is
in $\mathcal{C}^{1}[-a,a]$.  We consider the Cauchy problem
\begin{eqnarray}
D_{y}u\left(x,y\right)+F\left(x,y,u\left(x,y\right),
u\left(x,y\right),D_{x}u\left(x,y\right)\right)= 0\mbox{,
}\left(x,y\right)\in\Omega\label{Equation}
\end{eqnarray}
\begin{eqnarray}
u\left(x,0\right)=f\left(x\right)\mbox{,
}x\in\left(-a,a\right)\label{ICon}
\end{eqnarray}
Denote by $T:\mathcal{C}^{1}\left(\Omega\right) \rightarrow
\mathcal{C}^{0}\left(\Omega\right)$ the nonlinear partial
differential operator given by
\begin{eqnarray}
\begin{array}{ll}
\forall & u\in\mathcal{C}^{1}\left(\Omega\right)\mbox{ :} \\
\forall & \left(x,y\right)\in\Omega\mbox{ :} \\
& Tu:\left(x,y\right)\mapsto
D_{y}u\left(x,y\right)+F\left(x,y,u\left(x,y\right),
u\left(x,y\right),D_{x}u\left(x,y\right)\right) \\
\end{array}\label{PDEOpClassical}
\end{eqnarray}
Note that the equation (\ref{Equation}) may have several classical
solutions.  Indeed, in the particular case when the operator $T$
is linear and homogeneous, there is at least one classical
solution to (\ref{Equation}) which is the function which is
everywhere equal to $0$. However, the presence of the initial
condition (\ref{ICon}) may rule out some or all of the possible
classical solutions. What is more, there is a well known
\textit{physical} interest in \textit{nonclassical} or
\textit{generalized} solutions to (\ref{Equation}) through
(\ref{ICon}).  Such solution may, for instance, model shocks waves
in fluids.  In this regard, it is convenient to extend the PDE
operator $T$ to $\mathcal{ML}^{1}_{0}\left(\Omega\right)$ through
\begin{eqnarray}
\begin{array}{ll}
\forall & u\in\mathcal{ML}^{1}_{0}\left(\Omega\right)\mbox{ :} \\
\forall & \left(x,y\right)\in\Omega\mbox{ :} \\
& \mathcal{T}u:\left(x,y\right)\mapsto \left(I\circ S\right)\left(\mathcal{D}_{y}u+F\left(\cdot,\mathcal{D}_{x}u,u\right)\right)\left(x,y\right) \\
\end{array}\nonumber
\end{eqnarray}
As mentioned, the solution method for the Cauchy problem
(\ref{Equation}) through (\ref{ICon}) uses exactly the same
techniques that apply to the free problem \cite{vdWalt6}. However,
in order to incorporate the additional condition (\ref{ICon}), we
must adapt the method slightly.  In this regard, we consider the
space
\begin{eqnarray}
\mathcal{ML}^{1}_{0,0}\left(\Omega\right)=\{u\in
\mathcal{ML}^{1}_{0}\left(\Omega\right)\mbox{ :
}u\left(\cdot,0\right) \in\mathcal{C}^{1}[-a,a]\}\label{ML00Space}
\end{eqnarray}
and the mapping
\begin{eqnarray}
\mathcal{T}_{0}:\mathcal{ML}^{1}_{0,0}\left(\Omega\right)\ni
u\mapsto \left(\mathcal{T}u,\mathcal{R}_{0}u\right)\in
\mathcal{ML}^{0}_{0}\left(\Omega\right)\times
\mathcal{C}^{1}[-a,a]\label{IConOp}
\end{eqnarray}
where
\begin{eqnarray}
\begin{array}{ll}
\forall & u\in\mathcal{ML}^{1}_{0,0}\left(\Omega\right)\mbox{ :} \\
\forall & x\in [-a,a]\mbox{ :} \\
& \mathcal{R}_{0}u:x\mapsto u\left(x,0\right) \\
\end{array}\nonumber
\end{eqnarray}
That is, $\mathcal{R}^{0}$ assigns to each $u\in\mathcal{ML}^{1}_
{0,0}\left(\Omega\right)$ its restriction to
$\{\left(x,y\right)\in\Omega\mbox{ : }y=0\}$.  This amounts to a
\textit{separation} of the initial value problem (\ref{ICon}) form
the problem of satisfying the PDE (\ref{Equation}).

For the sake of a more compact exposition, we will denote by $X$
the space $\mathcal{ML}^{1}_ {0,0}\left(\Omega\right)$ and by $Y$
the space $\mathcal{ML}^{0}_{0}\left(\Omega\right)
\times\mathcal{C}^{1}[-a,a]$. On $\mathcal{C}^{1}[-a,a]$ we
consider the convergence structure $\lambda_{0}$, and with it the
associated u.c.s. $\mathcal{J}_{0}$.
\begin{definition}
For any $f\in\mathcal{C}^{1}[-a,a]$, and any filter $\mathcal{F}$
on $\mathcal{C}^{1}[-a,a]$,
\begin{eqnarray}
\mathcal{F}\in\lambda_{0}\left(f\right) \Leftrightarrow
[f]\subseteq \mathcal{F} \nonumber
\end{eqnarray}
Here $[x]$ denotes the filter generated by $x$.  That is,
\begin{eqnarray}
[x]=\{F\subseteq\mathcal{C}^{1}[-a,a]\mbox{ : }x\in F\}\nonumber
\end{eqnarray}
\end{definition}
The associated u.c.s. $\mathcal{J}_{0}$ on $\mathcal{C}^{1}[-a,a]$
consists of all filters $\mathcal{U}$ on
$\mathcal{C}^{1}[-a,a]\times \mathcal{C}^{1}[-a,a]$ that satisfies
\begin{eqnarray}
\begin{array}{ll}
\exists & k\in\mathbb{N}\mbox{ :} \\
& \left(\begin{array}{ll}
\forall & i=1,...,k\mbox{ :} \\
\exists & f_{i}\in\mathcal{C}^{1}[-a,a]\mbox{ :} \\
\exists & \mathcal{F}_{i}\mbox{ a filter on $\mathcal{C}^{1}[-a,a]$ :} \\
& \begin{array}{ll}
a) & \mathcal{F}_{i}\in\lambda_{0}\left(u_{i}\right) \\
b) & \left(\mathcal{F}_{1}\times\mathcal{F}_{1}\right)\cap...\cap\left(\mathcal{F}_{k}\times\mathcal{F}_{k}\right)\subseteq\mathcal{U} \\
\end{array} \\
\end{array}\right) \\
\end{array}
\end{eqnarray}
This u.c.s. is uniformly Hausdorff and complete.  The space $Y$
carries the product u.c.s. $\mathcal{J}_{Y}$ with respect to the
u.c.s's $\mathcal{J}_{a}$ on
$\mathcal{ML}^{0}_{0}\left(\Omega\right)$ and $\mathcal{J}_{0}$ on
$\mathcal{C}^{1}[-a,a]$.  That is, for any filter $\mathcal{V}$ on
$Y\times Y$
\begin{eqnarray}
\mathcal{V}\in\mathcal{J}_{Y}\Leftrightarrow \left(
\begin{array}{ll}
a) & \left(\pi_{0}\times\pi_{0}\right)\left(\mathcal{V}\right)\in\mathcal{J}_{a} \\
b) & \left(\pi_{1}\times\pi_{1}\right)\left(\mathcal{V}\right)\in\mathcal{J}_{0} \\
\end{array}\right)
\end{eqnarray}
Here $\pi_{0}$ denotes the projection on
$\mathcal{ML}^{0}_{0}\left(\Omega\right)$, and $\pi_{1}$ is the
projection on $\mathcal{C}^{1}[-a,a]$.  With this u.c.s., the
space $Y$ is uniformly Hausdorff and complete \cite[Proposition
2.3.3 (iii)]{Beattie}.

On the space $X$ we introduce an equivalence relation
$\sim_{\mathcal{T}_{0}}$ through
\begin{eqnarray}
\begin{array}{ll}
\forall & u,v\in X\mbox{ :} \\
& u\sim_{\mathcal{T}_{0}}v\Leftrightarrow \mathcal{T}_{0}u=\mathcal{T}_{0}v \\
\end{array}\label{T0Eq}
\end{eqnarray}
The quotient space $X/\sim_{\mathcal{T}_{0}}$ is denoted
$X_{\mathcal{T}_{0}}$.  With the mapping
$\mathcal{T}_{0}:X\rightarrow Y$ we may now associate an injective
mapping $\widehat{\mathcal{T}}_{0}:X_{\mathcal{T}_{0}}\rightarrow
Y$ so that the diagram\\ \\
\begin{math}
\setlength{\unitlength}{1cm} \thicklines
\begin{picture}(13,6)

\put(1.9,5.4){$X$} \put(2.3,5.5){\vector(1,0){6.0}}
\put(8.5,5.4){$Y$} \put(4.9,5.7){$\mathcal{T}_{0}$}
\put(2.0,5.2){\vector(0,-1){3.5}} \put(2.4,1.4){\vector(1,0){6.0}}
\put(1.9,1.3){$X_{\mathcal{T}_{0}}$} \put(8.5,1.3){$Y$}
\put(1.4,3.4){$q_{\mathcal{T}_{0}}$} \put(8.8,3.4){$i_{Y}$}
\put(8.6,5.2){\vector(0,-1){3.5}}
\put(4.9,1.6){$\widehat{\mathcal{T}}_{0}$}

\end{picture}
\end{math}\\
commutes.  Here $q_{\mathcal{T}_{0}}$ denotes the quotient mapping
associated with the equivalence relation $\sim_{\mathcal{T}_{0}}$,
and $i_{Y}$ is the identity mapping on $Y$.  We now define a
u.c.s. $\mathcal{J}_{\mathcal{T}_{0}}$ on $X_{\mathcal{T}_{0}}$ as
the \textit{initial} u.c.s. \cite[Proposition 2.1.1 (i)]{Beattie}
on $X_{\mathcal{T}_{0}}$ with respect to the mapping
$\widehat{\mathcal{T}}_{0}:X_{\mathcal{T}_{0}}\rightarrow Y$. That
is, \begin{eqnarray}
\begin{array}{ll}
\forall & \mathcal{U}\mbox{ a filter on $X_{\mathcal{T}_{0}}$ :} \\
& \mathcal{U}\in\mathcal{J}_{\mathcal{T}_{0}}\Leftrightarrow
\left(\widehat{\mathcal{T}}_{0}\times \widehat{\mathcal{T}}_{0}\right)\left(\mathcal{U}\right)\in\mathcal{J}_{Y} \\
\end{array}
\end{eqnarray}
Since $\widehat{\mathcal{T}}_{0}$ is injective, the u.c.s.
$\mathcal{J}_{\mathcal{T}_{0}}$ is uniformly Hausdorff, and
$\widehat{\mathcal{T}}_{0}$ is actually a uniformly continuous
embedding.  Moreover, if $X_{\mathcal{T}_{0}}^{\sharp}$ denotes
the completion of $X_{\mathcal{T}_{0}}$, then there exists a
unique uniformly continuous embedding
\begin{eqnarray}
\widehat{\mathcal{T}}_{0}^{\sharp}: X_{\mathcal{T}_{0}}^{\sharp}
\rightarrow Y\nonumber
\end{eqnarray}
such that the diagram\\ \\
\begin{math}
\setlength{\unitlength}{1cm} \thicklines
\begin{picture}(13,6)

\put(2.2,5.4){$X_{\mathcal{T}_{0}}$}
\put(2.8,5.5){\vector(1,0){5.8}} \put(8.8,5.4){$Y$}
\put(5.2,5.7){$\widehat{\mathcal{T}}_{0}$}
\put(2.5,5.2){\vector(0,-1){3.5}} \put(2.7,1.4){\vector(1,0){5.9}}
\put(2.1,1.2){$X_{\mathcal{T}_{0}}^{\sharp}$} \put(8.8,1.3){$Y$}
\put(1.6,3.4){$\iota_{X_{\mathcal{T}_{0}}}$}
\put(9.1,3.4){$i_{Y}$} \put(8.9,5.2){\vector(0,-1){3.5}}
\put(5.2,1.6){$\widehat{\mathcal{T}}_{0}^{\sharp}$}

\end{picture}
\end{math}\\
commutes.  Here $\iota_{X_{\mathcal{T}_{0}}}$ denotes the
uniformly continuous embedding associated with the completion
$X_{\mathcal{T}_{0}}^{\sharp}$ of $X_{\mathcal{T}_{0}}$. A
generalized solution to (\ref{Equation}) through (\ref{ICon}) is a
solution to the equation
\begin{eqnarray}
\widehat{\mathcal{T}}_{0}U^{\sharp}=\left(0,f\right)\label{GenEq}
\end{eqnarray}
The existence of generalized solutions is based on the following
basic approximation result \cite[Section 8]{Obergugenberger and
Rosinger}.
\begin{lemma}\label{Approx}
We have
\begin{eqnarray}
\begin{array}{ll}
\forall & \epsilon>0\mbox{ :} \\
\exists & \delta>0\mbox{ :} \\
\forall & \left(x_{0},y_{0}\right)\in\Omega\mbox{ :} \\
\exists & u=u_{\epsilon,x_{0},y_{0}}\in\mathcal{C}^{1}\left(\overline{\Omega}\right)\mbox{ :} \\
\forall & \left(x,y\right)\mbox{ :} \\
& \left(\begin{array}{l}
|x-x_{0}|<\delta \\
|y-y_{0}|<\delta \\
\end{array}\right)\Rightarrow -\epsilon\leq Tu\left(x,y\right)\leq 0 \\
\end{array}\label{1stApprox}
\end{eqnarray}
Furthermore, we also have
\begin{eqnarray}
\begin{array}{ll}
\forall & \epsilon>0\mbox{ :} \\
\exists & \delta>0\mbox{ :} \\
\forall & x_{0}\in[-a,a]\mbox{ :} \\
\exists & u=u_{\epsilon,x_{0}}\in\mathcal{C}^{1}\left(\overline{\Omega}\right)\mbox{ :} \\
& \begin{array}{ll}
a) & \forall \mbox{ }\mbox{ } \left(x,y\right)\in\overline{\Omega}\mbox{ :} \\
& \mbox{ }\mbox{ }\mbox{ }\left(\begin{array}{l}
|x-x_{0}|<\delta \\
|y|<\delta \\
\end{array}\right) \Rightarrow -\epsilon\leq Tu\left(x,y\right)\leq 0 \\
\end{array} \\
& \begin{array}{ll}
b) & u\left(x,0\right)=f\left(x\right)\mbox{, }x\in [-a,a]\mbox{, }|x-x_{0}|<\delta \\
\end{array}\\
\end{array}\label{2ndApprox}
\end{eqnarray}
\end{lemma}
As a consequence of the approximation result above, we now obtain
the \textit{existence} and \textit{uniqueness} of generalized
solutions to (\ref{Equation}) through (\ref{ICon}).  In this
regard, we introduce the concept of a finite initial adaptive
$\delta$-tiling.  A finite initial adaptive $\delta$-tiling of
$\Omega$ is any \textit{finite} collection
$\mathcal{K}_{\delta}=\{K_{1},...,K_{\nu}\}$ of perfect, compact
subsets of $\mathbb{R}^{2}$ with pairwise disjoint interiors such
that
\begin{eqnarray}
\begin{array}{ll}
\forall & K_{i}\in\mathcal{K}_{\delta}\mbox{ :} \\
& \left(x,y\right),\left(x_{0},y_{0}\right)\in K_{i}\Rightarrow \left(\begin{array}{l}
|x-x_{0}|<\delta \\
|y-y_{0}|<\delta \\
\end{array}\right) \\
\end{array}\label{FIADTiling1}
\end{eqnarray}
and
\begin{eqnarray}
\{\left(x,0\right)\mbox{ : }-a\leq x\leq a\}\cap\left(
\cup_{K_{i}\in\mathcal{K}_{\delta}}\partial K_{i}\right)\mbox{ at
most finite}\label{FIADTiling2}
\end{eqnarray}
where $\partial K_{i}$ denotes the boundary of $K_{i}$. For any
$\delta>0$ there exists at least one finite initial adaptive
$\delta$-tiling of $\Omega$, see for instance \cite[Section
8]{Obergugenberger and Rosinger}.
\begin{theorem}\label{ExTheorem}
For any $f\in\mathcal{C}^{1}[-a,a]$, there exists a unique
$U^{\sharp}\in X_{\mathcal{T}_{0}}^{\sharp}$ that satisfies
(\ref{GenEq}).
\end{theorem}
\begin{proof}
For every $n\in\mathbb{N}$, set $\epsilon_{n}=1/n$.  Applying
Lemma \ref{Approx}, we find $\delta_{n}>0$ such that
\begin{eqnarray}
\begin{array}{ll}
\forall & \left(x_{0},y_{0}\right)\in\Omega\mbox{ :} \\
\exists & u=u_{n,x_{0},y_{0}}\in\mathcal{C}^{1}\left(\overline{\Omega}\right)\mbox{ :} \\
\forall & \left(x,y\right)\mbox{ :} \\
& \left(\begin{array}{l}
|x-x_{0}|<\delta_{n} \\
|y-y_{0}|<\delta_{n} \\
\end{array}\right)\Rightarrow -\frac{\epsilon_{n}}{2}\leq Tu\left(x,y\right)\leq 0 \\
\end{array}\label{1stApprox}
\end{eqnarray}
and
\begin{eqnarray}
\begin{array}{ll}
\forall & x_{0}\in[-a,a]\mbox{ :} \\
\exists & u=u_{n,x_{0}}\in\mathcal{C}^{1}\left(\overline{\Omega}\right)\mbox{ :} \\
& \begin{array}{ll}
a) & \forall \mbox{ }\mbox{ } \left(x,y\right)\in\overline{\Omega}\mbox{ :} \\
& \mbox{ }\mbox{ }\mbox{ }\left(\begin{array}{l}
|x-x_{0}|<\delta \\
|y|<\delta \\
\end{array}\right) \Rightarrow -\frac{\epsilon_{n}}{2}\leq Tu\left(x,y\right)\leq 0 \\
\end{array} \\
& \begin{array}{ll}
b) & u\left(x,0\right)=f\left(x\right)\mbox{, }x\in [-a,a]\mbox{, }|x-x_{0}|<\delta \\
\end{array}\\
\end{array}\label{2ndApprox}
\end{eqnarray}
Let $\mathcal{K}_{\delta_{n}}=\{K_{1},...,K_{\nu_{n}}\}$ be a
finite initial adaptive $\delta_{n}$-tiling.  If $K_{i}\in
\mathcal{K}_{\delta_{n}}$, and
\begin{eqnarray}
K_{i}\cap\{\left(x,0\right)\mbox{ : }|x|\leq a\}=\emptyset
\end{eqnarray}
then take any
$\left(x_{0},y_{0}\right)\in\textrm{int}\left(K_{i}\right)$ and
set
\begin{eqnarray}
u^{i}_{n}=u_{n,x_{0},y_{0}}\nonumber
\end{eqnarray}
Otherwise, select $\left(x_{0},0\right)\in \left([-a,a]\times
\{0\}\right)\cap K_{i}$ and set
\begin{eqnarray}
u^{i}_{n}=u_{n,x_{0}} \nonumber
\end{eqnarray}
Consider the function $u_{n}\in\mathcal{ML}^{1}_{0,0}
\left(\Omega\right)$ defined as
\begin{eqnarray}
u_{n}=\left(I\circ
S\right)\left(\sum_{i=1}^{\nu}\chi_{i}u_{n}^{i}\right)\nonumber
\end{eqnarray}
where, for each $i$, $\chi_{i}$ is the indicator function of
$\textrm{int}\left(K_{i}\right)$.  It is clear that
\begin{eqnarray}
\begin{array}{ll}
\forall & \left(x,y\right)\in\Omega\mbox{ :} \\
& -\epsilon_{n}< \mathcal{T}u_{n}\left(x,y\right)\leq 0 \\
\end{array}\nonumber
\end{eqnarray}
and
\begin{eqnarray}
\begin{array}{ll}
\forall & x\in [-a,a]\mbox{ :} \\
& \mathcal{R}_{0}u_{n}\left(x\right)=0
\end{array}\nonumber
\end{eqnarray}
Let $U_{n}$ denote the $\sim_{\mathcal{T}_{0}}$ equivalence class
associated with $u_{n}$.  Then the sequence
$\left(\widehat{\mathcal{T}}_{0}U_{n}\right)$ converges to
$\left(0,f\right)\in Y$.  Since $\widehat{\mathcal{T}}_{0}$ is
uniformly continuous embedding, the sequence $\left(U_{n}\right)$
is a Cauchy sequence in $X_{\mathcal{T}_{0}}$ so that there exists
$U^{\sharp}\in X_{\mathcal{T}_{0}}^{\sharp}$ that satisfies
(\ref{GenEq}).  Moreover, this solution is unique since
$\widehat{\mathcal{T}}_{0}^{\sharp}$ is injective.
\end{proof}

\section{The Meaning of Generalized Solutions}

Regarding the meaning of the \textit{existence} and
\textit{uniqueness} of the generalized solution $U^{\sharp}\in
X_{\mathcal{T}_{0}}^{\sharp}$ to (\ref{Equation}) through
(\ref{ICon}), we recall the abstract construction of the
completion of a uniform convergence space.  Let
$\left(Z,\mathcal{J}\right)$ be a Hausdorff uniform convergence
space.  A filter $\mathcal{F}$ on $Z$ is a $\mathcal{J}$-Cauchy
filter whenever
\begin{eqnarray}
\mathcal{F}\times\mathcal{F}\in\mathcal{J}\nonumber
\end{eqnarray}
If $\textbf{C}\left(Z\right)$ denotes the collection of all
$\mathcal{J}$-Cauchy filters on $Z$, one may introduce an
equivalence relation $\sim_{c}$ on $\textbf{C}\left(Z\right)$
through
\begin{eqnarray}
\begin{array}{ll}
\forall & \mathcal{F},\mathcal{G}\in\textbf{C}\left(Z\right)\mbox{ :} \\
& \mathcal{F}\sim_{c}\mathcal{G}\Leftrightarrow \left(\begin{array}{ll}
\exists & \mathcal{H}\in\textbf{C}\left(Z\right)\mbox{ :} \\
& \mathcal{H}\subseteq \mathcal{G}\cap\mathcal{F} \\
\end{array}\right) \\
\end{array}
\end{eqnarray}
The quotient space $Z^{\sharp}=\textbf{C}\left(Z\right)/\sim_{c}$
serves as the underlying set of the completion of $Z$.  Note that,
since $Z$ is Hausdorff, if the filters
$\mathcal{F},\mathcal{G}\in\textbf{C}\left(Z\right)$ converge to
$x\in Z$ and $z\in Z$, respectively, then
\begin{eqnarray}
\mathcal{F}\sim_{c}\mathcal{G}\Leftrightarrow x=z\nonumber
\end{eqnarray}
so that one obtains an injective mapping
\begin{eqnarray}
\iota_{Z}:Z\ni z\mapsto [\lambda\left(z\right)]\in
Z^{\sharp}\nonumber
\end{eqnarray}
where $[\lambda\left(z\right)]$ denotes the equivalence class
generated by the filters which converge to $z\in Z$.  One may now
equip $Z^{\sharp}$ with a u.c.s. $\mathcal{J}^{\sharp}$ in such a
way that the mapping $\iota_{Z}$ is a uniformly continuous
embedding, $Z^{\sharp}$ is complete, and $\iota_{Z}\left(Z\right)$
is dense in $Z^{\sharp}$.

In the context of PDEs, and in particular the first order Cauchy
problem (\ref{Equation}) through (\ref{ICon}), the generalized
solution $U^{\sharp}\in X_{\mathcal{T}_{0}}$ to (\ref{Equation})
through (\ref{ICon}) may be expressed as
\begin{eqnarray}
U^{\sharp}=\left\{\mathcal{F}\in\textbf{C}\left(X_{\mathcal{T}_{0}}\right)\begin{array}{|ll}
a) & \pi_{0}\left(\widehat{\mathcal{T}}_{0}\left(\mathcal{F}\right)\right)\in\lambda_{a}\left(0\right) \\
& \\
b) & \pi_{1}\left(\widehat{\mathcal{T}}_{0}\left(\mathcal{F}\right)\right)\in\lambda_{0}\left(f\right) \\
\end{array}\right\}
\end{eqnarray}
Any classical solution $u\in\mathcal{C}^{1}\left(\Omega\right)$,
and more generally any shockwave solution
$u\in\mathcal{ML}^{1}_{0,0}\left(\Omega\right)=X$  to
(\ref{Equation}) through (\ref{ICon}), is mapped to the
generalized solution $U^{\sharp}$, as may be seen form the
following commutative diagram.\\ \\
\begin{math}
\setlength{\unitlength}{1cm} \thicklines
\begin{picture}(13,6)

\put(2.4,5.4){$X$} \put(2.8,5.5){\vector(1,0){5.8}}
\put(8.8,5.4){$Y$} \put(5.2,5.7){$\mathcal{T}_{0}$}
\put(2.5,5.2){\vector(0,-1){3.5}} \put(2.7,1.4){\vector(1,0){5.9}}
\put(2.1,1.2){$X_{\mathcal{T}_{0}}$} \put(8.8,1.3){$Y$}
\put(1.6,3.4){$q_{\mathcal{T}_{0}}$} \put(9.1,3.4){$i_{Y}$}
\put(8.9,5.2){\vector(0,-1){3.5}}
\put(5.2,1.6){$\widehat{\mathcal{T}}_{0}$} \put(2.4,
0.9){\vector(0,-1){3.5}} \put(8.9, 1.0){\vector(0,-1){3.5}}
\put(2.1,-3.0){$X_{\mathcal{T}_{0}}^{\sharp}$} \put(8.8,-3.0){$Y$}
\put(2.7,-2.9){\vector(1,0){5.9}}
\put(5.2,-2.7){$\widehat{\mathcal{T}}_{0}^{\sharp}$}
\put(1.6,-1.0){$\iota_{X_{\mathcal{T}_{0}}}$}
\put(9.1,-1.0){$i_{Y}$}
\end{picture}
\end{math}\\ \\ \\ \\ \\ \\ \\ \\

Hence there is a \textit{consistency} between the generalized
solutions $U^{\sharp}\in X_{\mathcal{T}_{0}}^{\sharp}$ and any
classical and shockwave solutions that may exists.

\section{Conclusion}

We have shown how the ideas developed in \cite{vdWalt5} through
\cite{vdWalt7} may be applied to solve initial and / or boundary
value problems for nonlinear systems of PDEs.  In this regard, we
have established the existence and uniqueness of generalized
solutions to a family of nonlinear first order Cauchy problems.
The generalized solutions are seen to be consistent with the usual
classical solutions, if the latter exists, as well as with shock
wave solutions. It should be noted that the above method applies
equally well to arbitrary  nonlinear \textit{systems} of
equations.

\end{document}